\documentclass{amsart}
\usepackage{amssymb} 
\usepackage{enumitem}
\usepackage{multirow}
\setlist[itemize]{leftmargin=12mm}
\setlist[enumerate]{leftmargin=12mm}
\usepackage{comment}

\DeclareMathOperator{\tr}{tr}
\DeclareMathOperator{\Aut}{Aut}
\DeclareMathOperator{\Cl}{Cl}

\DeclareMathOperator{\GL}{GL}
\DeclareMathOperator{\PGL}{PGL}

\DeclareMathOperator{\Gal}{Gal}
\DeclareMathOperator{\norm}{Norm}
\DeclareMathOperator{\ord}{\upsilon}

\DeclareMathOperator{\deltarel}{\delta_{\mathrm{rel}}}
\DeclareMathOperator{\deltasup}{\delta_{\mathrm{sup}}}
\newcommand{\Q}{{\mathbb Q}}
\newcommand{\Z}{{\mathbb Z}}

\newcommand{\F}{{\mathbb F}}
\newcommand{\PP}{{\mathbb P}}
\newcommand{\Nsf}{\mathbb{N}^\mathrm{sf}}
\newcommand{\cA}{\mathcal{A}}
\newcommand{\cU}{\mathcal{U}}
\newcommand{\cD}{\mathcal{D}}

\newcommand{\cG}{\mathcal{G}}

\newcommand{\cH}{\mathcal{H}}

\newcommand{\cL}{\mathcal{L}}
\newcommand{\cM}{\mathcal{M}}
\newcommand{\cN}{\mathcal{N}}

\newcommand{\sC}{\mathcal{C}}

\newcommand{\OO}{{\mathcal O}}
\newcommand{\ga}{\mathfrak{a}}

\newcommand{\ff}{\mathfrak{f}}
\newcommand{\gf}{\mathfrak{g}}
\newcommand{\gN}{\mathfrak{N}}
\newcommand{\fp}{\mathfrak{m}}
\newcommand{\fq}{\mathfrak{q}}
\newcommand{\mP}{\mathfrak{P}}     
\newcommand{\sS}{\mathfrak{S}}
\newcommand{\bs}{\mathbf{s}}

\begin {document}

\newtheorem{thm}{Theorem}
\newtheorem{lem}{Lemma}[section]

\newtheorem{cor}[lem]{Corollary}
\newtheorem{conj}{Conjecture}

\theoremstyle{definition}

\theoremstyle{remark}

\title[Asymptotic Fermat's Last Theorem]{
The Asymptotic 
Fermat's Last Theorem\\ for 
five-sixths of Real Quadratic Fields
}
\author{Nuno Freitas and Samir Siksek}
\address{
 Mathematisches Institut\\
Universit\"{a}t Bayreuth\\
95440 Bayreuth, Germany
}

\address{Mathematics Institute\\
	University of Warwick\\
Coventry\\
	CV4 7AL \\
	United Kingdom}


\date{\today}
\thanks{
The first-named author is supported
through a grant within the framework of the DFG Priority Programme 1489
{\em Algorithmic and Experimental Methods in Algebra, Geometry and Number Theory}.
The second-named
author is supported by an EPSRC Leadership Fellowship EP/G007268/1,
and EPSRC {\em LMF: L-Functions and Modular Forms} Programme Grant
EP/K034383/1.
}


\begin{abstract}
Let $K$ be a totally real field. 
By the {\em asymptotic Fermat's Last Theorem over $K$}
we mean the statement that there is a constant $B_K$
such that for prime exponent $p>B_K$ the only
 solutions to the Fermat equation 
\[
a^p+b^p+c^p=0, \qquad a,b,c \in K
\] 
are the trivial ones satisfying $abc=0$.
With the help of modularity, level lowering and
image of inertia comparisons we give an algorithmically
testable criterion which if satisfied by $K$ implies
the asymptotic Fermat's Last Theorem over $K$. 
Using techniques from analytic number theory, we show that our criterion
is satisfied by $K=\Q(\sqrt{d})$ for a subset
of $d \ge 2$ having density $5/6$ among the squarefree
positive integers. 
We can improve this density to $1$ if we assume
a standard
\lq\lq Eichler--Shimura\rq\rq\ conjecture. 
\end{abstract}


\maketitle



\section{Introduction}

\subsection{Historical background} 
Interest in the Fermat equation over various number
fields goes back to the 19th and early 20th Century.
For example, Dickson's {\em History of the Theory of Numbers}
\cite[pages 758 and 768]{Dickson}
mentions extensions by Maillet (1897) and Furtw\"{a}ngler (1910)
of classical ideas of Kummer to the Fermat equation
$x^p+y^p=z^p$ ($p>3$ prime) over the cyclotomic field
$\Q(\zeta_p)$. The ultimate work in this direction
 is due to Kolyvagin \cite{Kolyvagin} who, among
many other results, shows that if $x^p+y^p=z^p$
has a first case solution in $\Q(\zeta_p)$ then
$p^2 \mid (q^p-q)$ for all prime $q \le 89$
(here \lq first case\rq\ means that $(1-\zeta_p) \nmid xyz$).
For $d \ne 0$, $1$ squarefree, Hao and Parry \cite{HP}  
use the Kummer approach to prove several results on the Fermat equation
with prime exponent $p$ over $\Q(\sqrt{d})$,
subject to the condition that $p$ does not
divide the class number of $\Q(\sqrt{d},\zeta_p)$.
(See also \cite{HP} for references
to early work on the Fermat equation with various small 
exponents over quadratic and cubic fields.)
Others (e.g. \cite{DK}, \cite{Faddeev}, 
\cite{Tzermias}) 
treat the 
Fermat equation with fixed exponent $p$ as a curve, and determine the points
of low degree subject to restrictions on the Mordell--Weil
group of its Jacobian. A beautiful example
of this is due to Gross and Rohrlich \cite[Theorem 5]{GR} who
for $p=3$, $5$, $7$, $11$ determine the solutions to
$x^p+y^p=z^p$ over all number fields
$K$ of degree $\le (p-1)/2$.

However, the elementary,
cyclotomic and Mordell--Weil approaches to the Fermat equation have
had limited success. Indeed, even over $\Q$, no combination of these
approaches is known to yield a proof of Fermat's Last Theorem for
infinitely many prime exponents $p$.
We therefore look to
Wiles' proof of Fermat's Last Theorem \cite{TW}, \cite{Wiles}
for inspiration. 
Modularity of elliptic curves plays
a crucial r\^{o}le in the proof,
yet our understanding of modularity for elliptic curves
over general number fields is still embryonic.   
In contrast, there is  
a powerful modularity theory
for elliptic curves over totally real fields,
and it is natural to
apply this to the study of the Fermat equation
over totally real fields. 
The only work in this direction that we are aware of
is due to Jarvis and Meekin \cite{JMee}
who showed that the
Fermat equation $x^n+y^n=z^n$ has no solutions 
$x$, $y$, $z \in \Q(\sqrt{2})$ with $xyz \ne 0$ and $n \ge 4$.



\subsection{Our results}
In this paper we combine modularity and level lowering
with image of inertia comparisons and techniques
from analytic number theory to prove several results concerning the Fermat
equation over totally real fields. 
Let $K$ be a totally real field, and let $\OO_K$
be its ring of integers. By the \textbf{Fermat equation with exponent $p$ over $K$} we mean the equation
\begin{equation}\label{eqn:Fermat}
a^p+b^p+c^p=0, \qquad a,b,c\in \OO_K.
\end{equation}
A solution $(a,b,c)$ is called \textbf{trivial}
if $abc=0$, otherwise \textbf{non-trivial}.
The \textbf{asymptotic Fermat's Last Theorem over $K$}
is the statement that there is some bound $B_K$ such that for prime $p>B_K$,
all solutions  to the Fermat equation~\eqref{eqn:Fermat} are trivial.
If $B_K$ is effectively computable,
we shall refer to this as the
\textbf{effective asymptotic Fermat's Last Theorem over $K$}.

\begin{thm} \label{thm:FermatQuad}
Let $d\ge 2$ be squarefree, satisfying one of the following conditions
\begin{itemize}
\item[(i)] $d \equiv 3 \pmod{8}$,
\item[(ii)] $d \equiv 6$ or $10 \pmod{16}$, 
\item[(iii)] $d \equiv 2 \pmod{16}$ and $d$ has some prime
divisor $q \equiv 5$ or $7 \pmod{8}$,
\item[(iv)] $d \equiv 14 \pmod{16}$ and $d$ has some prime
divisor $q \equiv 3$ or $5 \pmod{8}$.
\end{itemize}
Then the effective asymptotic Fermat's Last Theorem holds over
$K=\Q(\sqrt{d})$.
\end{thm}

To state our other theorems, we need the following
standard conjecture, which of course 
is a generalization 
of the Eichler--Shimura Theorem
over $\Q$.
\begin{conj}[\lq\lq Eichler--Shimura\rq\rq]\label{conj:ES}
Let $K$ be a totally real field. Let $\ff$ be a Hilbert newform
of level $\cN$ and parallel weight $2$, and 
rational eigenvalues.
Then there is an elliptic curve $E_\ff/K$ with conductor $\cN$
having the same $\mathrm{L}$-function as $\ff$.
\end{conj}
\begin{thm}\label{thm:d5mod8}
Let $d>5$ be squarefree, satisfying $d \equiv 5 \pmod{8}$.
Write $K=\Q(\sqrt{d})$ and assume Conjecture~\ref{conj:ES} 
for $K$.
Then the effective asymptotic Fermat's Last Theorem holds over $K$.
\end{thm}

To state our general theorem, let $K$
be a totally real field, and let 
\begin{equation}\label{eqn:ST} 
\begin{gathered}
S=\{ \mP \; :\; \text{$\mP$ is a prime of $K$ above $2$}\}, \\
T=\{ \mP \in S \; : \; 
f(\mP/2)=1\},
\qquad 
U=\{ \mP \in S \; : \; 
3 \nmid \ord_\mP(2) \}.
\end{gathered}
\end{equation}
Here $f(\mP/2)$ denotes the residual degree of $\mP$.
We need an assumption, which we refer to throughout
the paper as (ES):
\[ \label{ES}
\text{\bf (ES)} \qquad
\left\{
\begin{array}{lll}
\text{either $[K:\Q]$ is odd;}\\
\text{or $T \ne \emptyset$;}\\
\text{or Conjecture~\ref{conj:ES} holds for $K$.}
\end{array}
\right.
\]

\begin{thm}\label{thm:FermatGen}
Let $K$ be a totally real field satisfying
 (ES). Let $S$, $T$ and $U$ be as in \eqref{eqn:ST}.
Write $\OO_S^*$ for the group of $S$-units of $K$.
Suppose that for every solution $(\lambda,\mu)$ to the $S$-unit equation 
\begin{equation}\label{eqn:sunit}
\lambda+\mu=1, \qquad \lambda,\, \mu \in \OO_S^* \, .
\end{equation}
there is 
\begin{enumerate}
\item[(A)] either some $\mP \in T$ that satisfies 
$\max\{ \lvert \ord_{\mP} (\lambda) \rvert, \lvert \ord_{\mP}(\mu) \rvert \} 
\le 4 \ord_{\mP}(2)$,
\item[(B)] or some $\mP \in U$ that satisfies both 
$\max\{ \lvert \ord_{\mP} (\lambda) \rvert, \lvert \ord_{\mP}(\mu) \rvert \} 
\le 4 \ord_{\mP}(2)$, and 
$\ord_{\mP}(\lambda \mu) \equiv \ord_{\mP}(2) \pmod{3}$.
\end{enumerate}
Then the asymptotic Fermat's Last Theorem holds over
$K$.
\end{thm}
We make the following remarks.
\begin{itemize}[leftmargin=25pt]
\item In contrast to Theorems~\ref{thm:FermatQuad} and~\ref{thm:d5mod8}, 
the constant 
$B_K$ implicit in the statement of Theorem~\ref{thm:FermatGen} 
is ineffective, though it is effectively computable if we
assume a suitable modularity statement, such as
modularity of elliptic curves over $K$ with full $2$-torsion.
Elliptic curves over real quadratic fields
are modular \cite{FHS},
and so 
the
implicit constants $B_K$ in Theorems~\ref{thm:FermatQuad},~\ref{thm:d5mod8} are
effectively computable.
\item By Siegel \cite{Siegel}, $S$-unit equations have 
a finite number of solutions.
These are effectively computable (e.g.
\cite{Smart}). Thus for any totally real $K$,
there is an algorithm for deciding whether the hypotheses of Theorem~\ref{thm:FermatGen}
are satisfied.
\item The $S$-unit equation \eqref{eqn:sunit} has precisely three 
solutions in $\Q \cap \OO_S^*$, 
namely $(\lambda,\mu)=(2,-1)$, $(-1,2)$, $(1/2,1/2)$.
We shall call these the \textbf{irrelevant solutions} to \eqref{eqn:sunit},
with other solutions being called \textbf{relevant}.
The irrelevant solutions
 satisfy (A) if $T \ne \emptyset$ and (B) if $U \ne \emptyset$.
\end{itemize}


It is natural to ask, for  fixed $n\ge 2$, for the
 \lq proportion\rq\
of totally real fields of degree $n$ 
such that the $S$-unit equation~\eqref{eqn:sunit}
has no relevant solutions. We answer this for $n=2$.
Specifically, let 
$\Nsf=\{ d \ge 2 \; : \; \text{$d$ squarefree} \}$.
The elements of $\Nsf$ are in bijection with the real
quadratic fields via $d \leftrightarrow K=\Q(\sqrt{d})$.
For a subset $\cU \subseteq \Nsf$, we shall define the
\textbf{relative density of $\cU$ in $\Nsf$} as
\begin{equation}\label{eqn:deltarel}
\deltarel(\cU)=\lim_{X \rightarrow \infty}
\frac{  \# \{ d \in \cU : d \le X \}  }{ \# \{ d \in \Nsf : d \le X\},}
\end{equation}
provided the limit exists. Let
\begin{equation}\label{eqn:cD}
\begin{gathered}
\sC=\{ d \in \Nsf :  \text{the $S$-unit equation
\eqref{eqn:sunit} has no relevant  solutions in $\Q(\sqrt{d})$} \} \\
\cD=\{ d \in \sC \; :\;  d \not \equiv 5 \pmod{8} \}.
\end{gathered}
\end{equation}

\begin{thm}\label{thm:density}
Let $\sC$ and $\cD$ be as above. Then
\begin{equation}\label{eqn:density}
\deltarel(\sC)=1, \qquad \deltarel(\cD)=5/6.
\end{equation}
If $d \in \cD$ then the effective
asymptotic Fermat's Last Theorem holds for $K=\Q(\sqrt{d})$.
The same
conclusion holds for $d \in \sC$ if we
assume Conjecture~\ref{conj:ES}.
\end{thm}
In other words, we are able to effectively bound the exponent in the Fermat
equation unconditionally for $5/6$ of real quadratic fields,
and assuming Conjecture~\ref{conj:ES} we can do this
for almost all real quadratic fields.


\subsection{Limitations of the original FLT strategy and variants}
We now answer the obvious question of how the proofs of the above  
differ from the proof of Fermat's Last Theorem over $\Q$. 
A basic sketch of the proof over $\Q$ is as follows. Suppose 
$a$, $b$, $c \in \Q$ satisfy $a^p+b^p+c^p=0$ 
with $abc \ne 0$ and $p \ge 5$ prime. 
Scale $a$, $b$, $c$ so that they become coprime integers.
After possibly permuting $a$, $b$, $c$ and 
changing signs we may suppose $a \equiv -1 \pmod{4}$
 and $2 \mid b$.
Consider the {\text Frey} elliptic curve
\begin{equation}\label{eqn:Frey}
E_{a,b,c} \; : \; Y^2=X(X-a^p)(X+b^p) 
\end{equation}
and 
let $\overline{\rho}_{E,p}$ be 
its mod $p$ Galois representation, where $E=E_{a,b,c}$.
Then $\overline{\rho}_{E,p}$ is irreducible
by Mazur \cite{Mazur} and  modular by Wiles and Taylor \cite{Wiles}, \cite{TW}.
Applying Ribet's level lowering theorem \cite{RibetLL}
shows that 
$\overline{\rho}_{E,p}$ 
arises from a weight $2$ newform of level $2$;
there are no such newforms 
giving a contradiction. 
Enough of modularity, irreducibility
and level lowering are now known for totally real fields that
we may attempt to carry out the same strategy 
over those fields. 
Let $K$ be a totally real field. 
If $a$, $b$, $c \in K$
satisfy $a^p+b^p+c^p=0$ 
and $abc \ne 0$, we can certainly scale $a$, $b$, $c$ so that
they belong to the ring of integers $\OO_K$, and thus $(a,b,c)$ is a non-trivial solution to \eqref{eqn:Fermat}. 
If the class number of $\OO_K$
is $>1$, we cannot assume coprimality
of $a$, $b$, $c$. 
We can however choose a finite set $\cH$
of prime ideals that represent the class group, and assume
(after suitable scaling)
that $a$, $b$, $c$ belong to $\OO_K$ and are coprime away from $\cH$.  
The Frey curve $E_{a,b,c}$ (again defined by \eqref{eqn:Frey}
but over $K$) is semistable outside $S \cup \cH$
(where $S$ is given in \eqref{eqn:ST}).
The non-triviality of the class group however obstructs
the construction of a Frey curve that is semistable outside $S$.
Applying level lowering to $\overline{\rho}_{E,p}$
yields a Hilbert newform $\ff$ of parallel weight $2$ 
and level divisible only by the primes in $S \cup \cH$.
In general, there are
newforms at these levels.
This situation is
analyzed by Jarvis and Meekin \cite{JMee} 
who find that
\begin{quote}
{\em \lq\lq\dots the numerology required to generalise the work of Ribet and
Wiles directly continues to hold for $\Q(\sqrt{2})$\dots
there are no other real quadratic fields for
which this is true \dots\rq\rq}
\end{quote}
It is helpful here to make a comparison with 
the equation $x^p+y^p+L^\alpha z^p=0$ over $\Q$,
with $L$ an odd prime, considered by Serre 
and Mazur \cite[p.\ 204]{SerreDuke}.
A non-trivial
solution to this latter equation gives rise, via modularity
and level lowering, to a classical weight $2$ newform $f$ of level
$2L$; for $L \ge 13$ there are such newforms 
and we face the same difficulty. Mazur however shows that
if $p$ is sufficiently large then $f$
corresponds to an elliptic curve $E^\prime$ with full $2$-torsion 
and conductor $2L$, and by classifying such elliptic curves
concludes that $L$ is either a Fermat or a Mersenne prime.
(Mazur's argument is unpublished, but can be found
in \cite[Section 15.5]{Cohen}.)
Mazur's argument adapted to our setting tells us that
a non-trivial solution to the Fermat equation over $K$
with $p$ sufficiently large gives rise to  
$E^\prime$ with full $2$-torsion and good reduction outside
$S \cup \cH$
(assumption (ES) is needed here).
It seems hopeless to classify
all such $E^\prime$ over all totally real fields
with a set of bad primes that varies with the  
class group
and might be arbitrarily large.

\subsection{Our approach over general totally real fields}
To go further we must somehow eliminate the bad primes
coming from the class group. Inspired by
Bennett and Skinner \cite{BenS}, and Kraus \cite{Kraus33p},
we study for $\fq \notin S$
the action of inertia groups $I_\fq$
on $E[p]$ for the Frey curve $E$. The curves $E$ and $E^\prime$
are related via 
$\overline{\rho}_{E,p} \sim \overline{\rho}_{E^\prime,p}$,
and so we obtain information about the
action of $I_\fq$ on $E^\prime[p]$. We use this to conclude
that $E^\prime$ has 
potentially good reduction
away from $S$. 

To summarize, assuming (ES),
a non-trivial solution to the Fermat equation
over $K$ with sufficiently large exponent $p$ yields an 
$E^\prime/K$ with full $2$-torsion and potentially
good reduction away from the set $S$ of primes above $2$. 
(There are such elliptic curves over every $K$, for example
the curve $Y^2=X^3-X$,
so we do not yet have a contradiction.)
However, such an $E^\prime$ can be represented by an $\OO_S$-point
on $Y(1)=X(1) \backslash \{\infty\}$ 
that pulls back to a $\OO_S$-point on $X(2) \backslash \{0,1,\infty\}$
(here $\OO_S$ is the ring of $S$-integers in $K$).
We can parametrize all such points (and therefore all such $E^\prime$)
 in terms of solutions
$(\lambda,\mu)$
to the $S$-unit equation \eqref{eqn:sunit}. 
This equation has solutions 
$(2,-1)$, $(-1,2)$, $(-1/2,-1/2)$ that we have called irrelevant above
(these correspond to $Y^2=X^3-X$),
and possibly others,
so we cannot yet conclude a contradiction.
However, 
the action of $I_\mP$ on $E[p]$,
for $\mP \in S$, 
gives
information on the valuations $\ord_\mP(\lambda)$, $\ord_\mP(\mu)$,
allowing us to prove Theorem~\ref{thm:FermatGen}.

\subsection{Our approach over real quadratic fields}
Next we specialize to real quadratic fields $K=\Q(\sqrt{d})$,
and we would like to understand solutions $(\lambda,\mu)$
to the $S$-unit equation \eqref{eqn:sunit}. By considering
$\norm(\lambda)$ and $\norm(\mu)$, which must both be of the form
$\pm 2^r$, 
we show that
relevant $(\lambda,\mu)$ give rise to solutions to the
equation
\begin{equation}\label{eqn:expdio}
(\eta_1 \cdot 2^{r_1} -\eta_2 \cdot 2^{r_2} +1)^2 
-\eta_1 \cdot 2^{r_1+2}= d v^2
\end{equation}
with $\eta_1=\pm 1$, $\eta_2=\pm 1$ and $r_1 \ge r_2 \ge 0$
and $v \in \Z\backslash \{0\}$. This approach has the merit of eliminating
the field $K$, since all the unknowns are now rational
integers. We solve this equation completely
for the $d$ appearing in Theorems~\ref{thm:FermatQuad}
and~\ref{thm:d5mod8};  there are a handful
of solutions and these satisfy conditions (A) or (B) of 
Theorem~\ref{thm:FermatGen}. This gives proofs for
Theorems~\ref{thm:FermatQuad} and~\ref{thm:d5mod8}.

To prove Theorem~\ref{thm:density} we need to show that
the set of squarefree $d \ge 2$ for which equation~\eqref{eqn:expdio}
has a solution is of density $0$. 
For this we employ sieving modulo Mersenne numbers
$M_m=2^m-1$. Modulo $M_m$ there are at most $4 m^2$ possibilities for
the left-hand side of \eqref{eqn:expdio}. The number of
possibilities for squares modulo $M_m$ is roughly at most
$M_m/2^{\omega(M_m)}$ where $\omega(n)$ denotes the number
of prime divisors of $n$. \textbf{If we can invert $v^2$ modulo $M_m$}
we will conclude that $d$ belongs to roughly at most $4 m^2 M_m/2^{\omega(M_m)}$
congruence classes modulo $M_m$ and so the set of possible
$d$ has density roughly at most $4m^2/2^{\omega(M_m)}$.
Using the {\em Prime Number Theorem} and the {\em Primitive Divisor Theorem}
it is possible to choose values
of $m$ so that this ratio tends to $0$.
 However there
is no reason to suppose that $v$ and $M_m$
are coprime and the argument
needs to be combined with other techniques from analytic
number theory 
to prove  
that the density of possible $d$ is
indeed $0$.

\subsection{Relation to other equations of Fermat-type and to modular curves}
We place our work in the wider
Diophantine context by mentioning related results on Fermat-type
equations. 
Many results of the FLT or asymptotic FLT kind (but over $\Q$)
are established by Kraus \cite{KrausFermat}, 
Bennett and Skinner \cite{BenS}, and Bennett, Vatsal and
Yazdani \cite{BVY} for equations of the form 
$A x^p+B y^p+C z^p=0$, $A x^p+B y^p=C z^2$, $A x^p+B y^p=C z^3$.
Perhaps the most satisfying work on Fermat-type equations over $\Q$ 
is the paper of Halberstadt and Kraus \cite{HK}, in which they
show, for any triple of  odd pairwise coprime integers $A$, $B$, $C$, 
there is a set
of primes $p$ of positive density, such that all solutions to
$A x^p+B y^p+C z^p=0$ are trivial. 

One can also consider the
analogy between the Fermat curves and various 
modular curves such as $X_0(p)$, $X_1(p)$. In
this framework Mazur's theorems \cite{Mazur}
 are analogues of Fermat's Last Theorem
over $\Q$. Merel's Uniform Boundedness Theorem \cite{Merel} states that for
a number field $K$, there is a bound $B_n$ that depends only
on the degree $n=[K:\Q]$, such that for $p>B_n$ the $K$-points
on $X_1(p)$ are cusps. Interestingly, our bounds $B_K$
for the asymptotic Fermat (when applicable and effective)
depend on the totally real field, and indeed the Hilbert newforms
at certain levels.
The links between Fermat-type curves and modular curves are 
far deeper than one might suspect. 
For example, to solve $x^p+y^p=2 z^p$,
$x^p+y^p=z^2$ and $x^p+y^p=z^3$ (and thereby complete the
resolution of $x^p+y^p=2^\alpha z^p$ that was started 
by Ribet \cite{RibetDio}),
Darmon and Merel \cite{DM} not only 
needed the theorems of Mazur, Ribet and Wiles,
but were also forced to study the rational points on the 
 modular curves $X_{\mathrm{ns}}^+(p) \times_{X(1)} X_0(r)$,
with $r=2$, $3$. 
In the reverse direction, methods of the present
paper are used in \cite{AS} to prove asymptotic results for semistable points on 
$X_{\mathrm{ns}}^+(p)$ and $X_{\mathrm{s}}^+(p)$ over
totally real fields.

\subsection{Notational conventions}\label{sec:conv}
Throughout $p$ denotes a rational prime, and $K$ a 
totally real number field, with 
ring of integers $\OO_K$.
For a non-zero ideal $\ga$ of $\OO_K$, we denote by
$[\ga]$ the class of $\ga$ in the class group $\Cl(K)$. 
For a non-trivial solution $(a,b,c)$ to the Fermat equation \eqref{eqn:Fermat},
let
\begin{equation}\label{eqn:cG}
\cG_{a,b,c}:=a \OO_K+b\OO_K+c\OO_K \, .
\end{equation}
and let $[a,b,c]$ denote
the class of $\cG_{a,b,c}$ in $\Cl(K)$.
We exploit the well-known fact
(e.g. \cite[Theorem VIII.4]{CassFr})
that every ideal class contains
infinitely many prime ideals.
Let $\mathfrak{c}_1,\dotsc,\mathfrak{c}_h$
be the ideal classes of $K$.
For each class $\mathfrak{c}_i$, we choose 
(and fix) a prime ideal $\fp_i \nmid 2$ of smallest possible
norm representing $\mathfrak{c}_i$. Throughout,
the set $\cH$ denotes our fixed choice of odd prime ideals
representing the class group:
$\cH=\{\fp_1,\dots,\fp_h\}$.
The sets $S$, $T$  and $U$ are given in \eqref{eqn:ST}. Observe
that $S \cap \cH=\emptyset$.

Let $G_K=\Gal(\overline{K}/K)$.
For an elliptic curve $E/K$,
we write
\begin{equation}\label{eqn:rho}
\overline{\rho}_{E,p}\; :\; G_K \rightarrow \Aut(E[p]) \cong \GL_2(\F_p),
\end{equation}
for the representation of $G_K$ on the $p$-torsion of $E$. 
For
a Hilbert eigenform $\ff$ over $K$, we let $\Q_\ff$ denote
the field generated by its eigenvalues.
In this situation $\varpi$ will denote a prime of $\Q_\ff$
above $p$; of course if $\Q_\ff=\Q$ we write
$p$ instead of $\varpi$. All other primes we consider are primes
of $K$. 
We reserve the symbol 
$\mP$ for  primes belonging to $S$, and 
$\fp$ for 
primes belonging to $\cH$.
An arbitrary prime of $K$ is denoted by $\fq$,  and $G_\fq$
and $I_\fq$ are the decomposition and inertia subgroups of $G_K$ 
at $\fq$.

\subsection*{Acknowledgements}
We are indebted to Alex Bartel, Frank Calegari, John Cremona, 
Lassina Demb\'{e}l\'{e},
Fred Diamond,
Tim Dokchitser, David Loeffler, Michael Stoll
 and Panagiotis Tsaknias
for useful discussions.

\section{Theoretical Background}
In this section we summarize the theoretical
results we need for modularity,
irreducibility of mod $p$ Galois representations, level lowering,
and Conjecture~\ref{conj:ES}.
\subsection{Modularity of the Frey Curve}
We shall need the following 
recent \cite{FHS} special case of the modularity conjecture
for elliptic curves over totally real fields.

\begin{thm}
\label{thm:modgen}
Let $K$ be a totally real field. Up to isomorphism over $\overline{K}$,
there are at most finitely many non-modular 
elliptic curves $E$ over $K$. Moreover, if $K$ is real quadratic,
then all elliptic curves over $K$ are modular.
\end{thm}



\begin{cor}\label{cor:Freymod}
Let $K$ be a totally real field.
There is some constant $A_K$ depending only on $K$, 
such that for any non-trivial solution 
$(a,b,c)$ of the Fermat equation \eqref{eqn:Fermat}
with prime exponent $p>A_K$,
 the Frey curve $E_{a,b,c}$ given by \eqref{eqn:Frey} is modular.
\end{cor}
\begin{proof}
By Theorem~\ref{thm:modgen}, there are at most finitely many possible
$\overline{K}$-isomorphism classes of elliptic curves over $K$
that are non-modular. Let $j_1,\dots,j_n \in K$ be the 
$j$-invariants of these classes. 
Write $\lambda=-b^p/a^p$. The $j$-invariant of $E_{a,b,c}$ is 
\[
j(\lambda)=2^8 \cdot (\lambda^2-\lambda+1)^3 \cdot \lambda^{-2} (\lambda-1)^{-2}.
\]
Each equation $j(\lambda)=j_i$ has at most six solutions $\lambda \in K$. 
Thus there are values $\lambda_1,\dots,\lambda_m \in K$
(with $m \le 6n$) such that if $\lambda \ne \lambda_k$ for 
all $k$ then $E_{a,b,c}$ is modular. 
If $\lambda=\lambda_k$ then 
\[
(-b/a)^p=\lambda_k, \qquad (-c/a)^p=1-\lambda_k. 
\]
This pair of equations results in a bound for $p$
unless $\lambda_k$ and $1-\lambda_k$ are both
roots of unity, which is impossible as $K$
is real, and so the only roots of unity are $\pm 1$. 
\end{proof}

\noindent \textbf{Remark.} The constant $A_K$ is 
ineffective: in \cite{FHS} 
it is shown that an elliptic curve $E$ over a totally real
field $K$ is modular except possibly if it gives rise to a $K$-point
on one of handful of modular curves of genus $\ge 2$,
and  
Faltings' Theorem \cite{Faltings1} (which is ineffective)
gives the finiteness. 
If 
$K$ is quadratic,
we can take $A_K=0$.

\subsection{Irreducibility of mod $p$ representations of elliptic curves}
We need the following,  derived in
\cite[Theorem 2]{FSirred} from the
work of David \cite{DavidI}
and Momose \cite{Momose}, who in turn build  
on Merel's Uniform Boundedness Theorem~\cite{Merel}.
\begin{thm}
\label{thm:irred}
Let $K$ be a Galois totally real field. 
There is an effective constant $C_K$, depending only on $K$, such that the 
following holds. 
If $p > C_K$ is prime,
and  $E$ is an elliptic curve over $K$ which is semistable at all $\fq \mid p$,
then $\overline{\rho}_{E,p}$ is irreducible.
\end{thm}

\subsection{Level Lowering}
We need a level lowering
result that plays the r\^ole of the Ribet step \cite{RibetLL}
in the proof of Fermat's Last Theorem. Fortunately such a result follows by 
 combining work of
Fujiwara \cite{Fuj},
Jarvis \cite{Jarv} and 
Rajaei \cite{Raj}. 

\begin{thm}[level lowering]\label{thm:levell} 
Let $K$ be a totally real field, and $E/K$ an elliptic curve of conductor $\cN$.
Let $p$ be a rational prime.
For a prime ideal $\fq$ of $K$ denote by $\Delta_\fq$ the discriminant of a local
minimal model for $E$ at $\fq$.
Let
\begin{equation}\label{eqn:Np}
\cM_p := \prod_{
\substack{\fq \Vert \cN,\\ p \mid \ord_\fq(\Delta_\fq)}
} {\fq}, \qquad\quad \cN_p:=\frac{\cN}{\cM_p} \, .
\end{equation}
Suppose the following
\begin{enumerate}
\item[(i)] $p\ge 5$,  
the ramification index $e(\fq/p) < p-1$ for all $\fq \mid p$, and
$\Q(\zeta_p)^{+} \not \subseteq K$.
\item[(ii)] $E$ is modular,
\item[(iii)] $\overline{\rho}_{E,p}$ is irreducible,
\item[(iv)] $E$ is semistable at all $\fq \mid p$,
\item[(v)]  $p \mid \ord_\fq(\Delta_\fq)$ for all $\fq \mid p$. 
\end{enumerate}
Then, there is a Hilbert eigenform $\ff$ 
of parallel weight $2$ that is new at level $\cN_p$ and some prime $\varpi$ of $\Q_\ff$
such that $\varpi \mid p$
and $\overline{\rho}_{E,p} \sim \overline{\rho}_{\ff,\varpi}$.
\end{thm}
\begin{proof} 
As noted above, we make use of  
the theorems in \cite{Fuj},
\cite{Jarv},
\cite{Raj}. Assumption (i) takes care of some 
technical restrictions in those theorems.

By assumption (ii), there is a newform $\ff_0$ of parallel weight $2$, level $\cN$ and
field of coefficients $\Q_{\ff_0} = \Q$, such that $\rho_{E,p} \sim \rho_{\ff_0,p}$.
Thus $\overline{\rho}_{E,p}$ is modular, and by (iii) is irreducible.
Since $K$ may
be of even degree, in order to apply the main result of \cite{Raj}, we need to
add an auxiliary (special or supercuspidal) prime to the level. From
\cite[Theorem 5]{Raj} we can add an auxiliary (special) prime
$\mathfrak{q}_0 \nmid \cN$ so that
$\overline{\rho}_{\ff_0,p}(\sigma_{\fq_0})$ is conjugate to
$\overline{\rho}_{\ff_0,p}(\sigma)$, where $\sigma_{\fq_0}$ denotes a Frobenius element of $G_K$
at $\fq_0$ and $\sigma$ is complex conjugation. We now
apply the main theorem of \cite{Raj} to remove from the level all primes $\fq
\nmid p$ dividing $\cM_p$. Next we  remove from the level the primes above
$p$ without changing the weight. 
By \cite[Theorem 6.2]{Jarv} we can do this provided
$\overline{\rho}_{E,p}\vert_{G_{\fq}}$ is finite at all $\fq \mid p$, 
where $G_{\fq}$
is the decomposition subgroup of $G_K$ at $\fq$. 
But from (iv), $\fq$ is a prime of good or multiplicative reduction for $E$.
In the former case,
$\overline{\rho}_{E,p}\vert_{G_{\fq}}$ is finite; in the latter case it is finite
by (v).
Finally, from the condition imposed on
$\mathfrak{q}_0$ it follows that $\norm(\mathfrak{q}_0) \not\equiv 1 \pmod{p}$, 
and we can apply Fujiwara's version of Mazur's principle \cite{Fuj} to remove
$\mathfrak{q}_0$ from the level. We conclude that there is an eigenform 
$\ff$ of parallel weight $2$, new at level $\cN_p$, and a prime $\varpi \mid p$ of $\Q_{\ff}$ such that
$\overline{\rho}_{E,p} \sim \overline{\rho}_{\ff_0,p} \sim \overline{\rho}_{\ff,\varpi}$. 
\end{proof}


\subsection{Eichler--Shimura}\label{sec:ES}
The following is partial result 
towards Conjecture~\ref{conj:ES}.
\begin{thm}[Blasius, Hida] 
\label{thm:Zhang}
Let $K$ be a totally real field and let $\ff$ be a Hilbert newform 
over $K$ of level $\cN$ and parallel weight $2$, such that
$\Q_\ff=\Q$. 
Suppose that 
\begin{enumerate}
\item[(a)] either $[K:\Q]$ is odd,
\item[(b)] or there is a finite prime $\fq$ 
such that $\pi_\fq$ belongs to the discrete series,
where $\pi$ is the cuspidal automorphic representation of
$\GL_2(\mathbb{A}_K)$ attached to $\ff$.
\end{enumerate}
Then there is an elliptic curve $E_\ff/K$ of conductor $\cN$
with the same $\mathrm{L}$-function as~$\ff$.
\end{thm}
Theorem~\ref{thm:Zhang} is derived by
Blasius \cite{Blasius}
from the work of Hida \cite{Hida}.
If $\ord_\fq(\cN)=1$ then  
$\pi_\fq$ is special and (b)
is satisfied. This is the only case of (b) that we use;
proofs of this case of Theorem~\ref{thm:Zhang}
are given by Darmon \cite{Darmon} and Zhang \cite{Zhang}.

\begin{cor}\label{cor:ES}
Let $E$ be an elliptic curve over a totally real field $K$
and $p$ an odd prime. 
Suppose $\overline{\rho}_{E,p}$ is irreducible,
and $\overline{\rho}_{E,p} \sim \overline{\rho}_{\ff,p}$
for some Hilbert newform $\ff$ over $K$ of parallel weight $2$ 
with $\Q_\ff=\Q$. 
Let $\fq \nmid p$ be a prime of $K$ such that
\begin{enumerate}
\item[(a)] $E$ has potentially multiplicative reduction at $\fq$;
\item[(b)] $p \mid \# \overline{\rho}_{E,p}(I_\fq)$;
\item[(c)] $p \nmid (\norm_{K/\Q}(\fq) \pm 1)$.
\end{enumerate} 
Then there is an elliptic curve $E_\ff/K$ of conductor $\cN$
with the same $\mathrm{L}$-function as~$\ff$.
\end{cor}
\begin{proof}
Write $c_4$ and $c_6$ for the usual $c$-invariants of $E$,
which are non-zero as
$E$ has potentially multiplicative reduction
at $\fq$. Let $\gamma=-c_4/c_6$.
Write $\chi$ for the quadratic character
associated to $K(\sqrt{\gamma})/K$ and $E\otimes \chi$ 
for the $\gamma$-quadratic twist of $E$. 
 By \cite[Theorem V.5.3]{SilvermanII},
$E\otimes \chi$ has split multiplicative reduction at $\fq$.
Let $\gf=\ff \otimes \chi$. 
As $\chi$ is quadratic and $\Q_\ff=\Q$ we have $\Q_\gf=\Q$. 

Suppose $\gf$ is new at level $\cN_\gf$.
We will prove 
$\ord_\fq(\cN_\gf)=1$. Then, by
Theorem~\ref{thm:Zhang} there is an elliptic curve
$E_\gf$ over $K$ having the same $\mathrm{L}$-function as $\gf$.
Thus the $\mathrm{L}$-functions of $E_\gf \otimes \chi$ and $\gf \otimes \chi=\ff$
are equal, and we take $E_\ff=E_\gf \otimes \chi$.

It remains to prove $\ord_\fq(\cN_\gf)=1$. Since
$\overline{\rho}_{E\otimes \chi,p} \sim \overline{\rho}_{\gf,p}$,
the two representations have the same optimal Serre level
$\gN$ (say). 
Now $E\otimes \chi$ has multiplicative reduction at $\fq$, so
$\ord_\fq(\gN)=0$ or $1$. 
Since $E$ and $E\otimes \chi$ are isomorphic over $K(\sqrt{\gamma})$,
and as $p \mid \# \overline{\rho}_{E,p}(I_\fq)$,
we have
$p \mid \# \overline{\rho}_{E\otimes \chi,p}(I_\fq)$.
Hence $\ord_{\fq}(\gN) \ne 0$, and so
$\ord_{\fq}(\gN)=1$.

We  now think of $\gN$ as the optimal Serre level $\overline{\rho}_{\gf, p}$
and compare it to the level $\cN_\gf$ of $\gf$. 
By \cite[Theorem 1.5]{Jarvis},
$\ord_\fq(\gN)=\ord_\fq(\cN_\gf)$, except possibly when $\ord_\fq(\cN_\gf)=1$
and $\ord_\fq(\gN)=0$ or when $\norm_{K/\Q}(\fq) \equiv \pm 1 \pmod{p}$.
The former is impossible as $\ord_\fq(\gN)=1$ and the latter is ruled out
by (c).
Thus $\ord_\fq(\cN_\gf)=1$.
\end{proof}

\section{Computations}
\subsection{Behaviour at odd primes} 
For $u$, $v$, $w \in \OO_K$ such that $uvw\ne 0$ and $u+v+w=0$,
let
\begin{equation}\label{eqn:2tors}
E: y^2=x(x-u)(x+v).
\end{equation}
The invariants $c_4$, $c_6$, $\Delta$, $j$ have their usual meanings 
and are given by:
\begin{equation}\label{eqn:inv}
\begin{gathered}
c_4=16(u^2-vw)=16(v^2-wu)=16(w^2-uv),\\
c_6=-32(u-v)(v-w)(w-u), \qquad 
\Delta=16 u^2 v^2 w^2, \qquad j=c_4^3/\Delta \, .
\end{gathered}
\end{equation}

The following elementary lemma is 
a straightforward consequence of the properties of elliptic
curves over local fields (e.g.\ \cite[Sections VII.1 and VII.5]{SilvermanI}).
\begin{lem}\label{lem:elem}
With the above notation, let $\fq \nmid 2$ be a prime and let
\[
s=\min\{\ord_\fq(u),\ord_\fq(v),\ord_\fq(w) \}.
\]
Write $E_{\mathrm{min}}$ for a local minimal model at $\fq$.
\begin{enumerate}
\item[(i)] $E_{\mathrm{min}}$ has good reduction at $\fq$ if and only if
$s$ is even and 
\begin{equation}\label{eqn:uvweq}
\ord_\fq(u)=\ord_\fq(v)=\ord_\fq(w).
\end{equation}
\item[(ii)] $E_{\mathrm{min}}$ has multiplicative reduction at $\fq$
if and only if $s$ is even and \eqref{eqn:uvweq} fails to hold. In this
case the minimal discriminant 
$\Delta_{\fq}$ at $\fq$ satisfies 
\[
\ord_\fq(\Delta_{\fq})=2 \ord_\fq(u)+2\ord_\fq(v)+2 \ord_\fq(w)-6s.
\]
\item[(iii)] $E_{\mathrm{min}}$ has
additive reduction if and only if $s$ is odd.
\end{enumerate}
\end{lem}

\subsection{Conductor of the Frey curve}
Let $(a,b,c)$ be a non-trivial solution to the 
Fermat equation~\eqref{eqn:Fermat}.
Let $\cG_{a,b,c}$ be as given in \eqref{eqn:cG},
which we think of as the greatest
common divisor of $a$, $b$, $c$. 
An odd prime not dividing $\cG_{a,b,c}$ is a prime
of good or multiplicative reduction for
$E_{a,b,c}$ and does not
appear in the final level $\cN_p$, as we see in due course.
An odd prime dividing
$\cG_{a,b,c}$ exactly once is an additive prime, and 
does appear in $\cN_p$.
To control $\cN_p$, we need to
control $\cG_{a,b,c}$. The following lemma achieves this. 
We refer to Section~\ref{sec:conv}
for the notation. 
\begin{lem}\label{lem:gcd}
Let $(a,b,c)$ be a non-trivial solution to \eqref{eqn:Fermat}.
There is a non-trivial integral solution $(a^\prime,b^\prime,c^\prime)$ to \eqref{eqn:Fermat} such 
that the following hold.
\begin{enumerate}
\item[(i)] For some $\xi \in K^*$, we have
$a^\prime= \xi a$, $b^\prime= \xi b$, $c^\prime=\xi c$.
\item[(ii)] $\cG_{a^\prime,b^\prime,c^\prime}=\fp$ for some $\fp \in \cH$.
\item[(iii)] $[a^\prime,b^\prime,c^\prime]=[a,b,c]$.
\end{enumerate}
\end{lem}
\begin{proof}
Let $\fp \in \cH$ satisfy $[\mathcal{G}_{a,b,c}]=[\fp]$,
so $\fp=(\xi) \cdot \mathcal{G}_{a,b,c}$ for some $\xi \in K^*$. 
Let $a^\prime$, $b^\prime$,
$c^\prime$ be as in (i).  
Note 
$(a^\prime)=(\xi) \cdot (a)=\fp \cdot \cG_{a,b,c}^{-1} (a)$
which is an integral ideal, since $\cG_{a,b,c}$ (by its
definition) divides $a$.
Thus $a^\prime$ is in $\OO_K$ and similarly so are $b^\prime$ and $c^\prime$.
For (ii) and (iii), note that
\[
\cG_{a^\prime,b^\prime,c^\prime}
=a^\prime \OO_K+b^\prime \OO_K+c^\prime \OO_K
=(\xi)\cdot (a \OO_K+b \OO_K+c \OO_K)
=(\xi)\cdot \cG_{a,b,c}=\fp.
\]
\end{proof}

\begin{lem}\label{lem:cond}
Let 
$(a,b,c)$ be a non-trivial solution to the 
Fermat equation \eqref{eqn:Fermat} with prime exponent $p$
satisfying $\cG_{a,b,c} = \fp$, where $\fp \in \cH$. 
Write $E$
for the Frey curve in \eqref{eqn:Frey}, and let $\Delta$ be its discriminant. 
Then at all $\fq \notin S \cup \{\fp\}$, 
the model $E$ is minimal, semistable, and satisfies 
$p \mid \ord_\fq(\Delta)$.
Let $\cN$ be the conductor of $E$, and let $\cN_p$
be as defined in \eqref{eqn:Np}. Then
\begin{equation}\label{eqn:cnp}
\cN=
\fp^{s_\fp}\cdot
\prod_{\mP \in S} \mP^{r_\mP}
\cdot 
\prod_{\substack{\fq \mid abc \\ 
\fq \notin S \cup \{\fp\} 
}} \fq, 
\qquad 
\qquad
\cN_p= 
\fp^{s_\fp^\prime}\cdot
\prod_{\mP \in S} 
\mP^{r_\mP^\prime}
 \, ,
\end{equation}
where $0 \le r_\mP^\prime \le r_\mP \le 
2 + 6\ord_{\mP} (2)$ and $0 \le s_\fp^\prime \le s_\fp \le 2$. 
\end{lem}
\begin{proof}
Suppose first that $\fq \mid abc$ and $\fq \notin S \cup \{\fp\}$.
Since $\cG_{a,b,c}= \fp$, the prime $\fq$ divides precisely one of $a$, $b$, $c$.
From \eqref{eqn:inv}, $\fq \nmid c_4$ so the model \eqref{eqn:Frey}
is minimal and
has multiplicative reduction at $\fq$, and $p \mid \ord_\fq(\Delta)$.
By \eqref{eqn:Np}, we see that $\fq \nmid \cN_p$.

For $\mP \in S$
we have $r_\mP=\ord_\mP(\cN) \le 2+6 \ord_{\mP}(2)$ by \cite[Theorem IV.10.4]{SilvermanII}.
We observe by \eqref{eqn:Np} that $r_\mP^\prime=r_\mP$
unless $E$ has multiplicative reduction at $\mP$ 
and $p \mid \ord_\mP(\Delta_\mP)$ in which case 
$r_\mP=1$ and $r_\mP^\prime=0$. 
Finally, recall by our choice of class group representatives $\cH$,
we have $\fp \nmid 2$.
As $E$ has full $2$-torsion over $K$, the wild part of the conductor
of $E/K$ at $\fp$ vanishes (see \cite[page 380]{SilvermanII}).
Thus $s_\fp=\ord_\fq(\cN) \le 2$. 
Again $s_\fp^\prime=s_\fp$ unless $E$ has multiplicative
reduction at $\fp$ and $p \mid \ord_\fp(\Delta_\fp)$
in which case $s_\fp=1$ and $s_\fp^\prime=0$.
\end{proof}

\subsection{Images of inertia}\label{sec:inertia}
We gather information needed regarding
images of inertia $\overline{\rho}_{E,p}(I_\fq)$. 
This is crucial for
applying Corollary~\ref{cor:ES}, and for controlling the 
behaviour at the primes in $S \cup \{ \fp\}$ of the newform
obtained by level lowering. 
\begin{lem}\label{lem:inertiaGeneral}
Let $E$ be an elliptic curve over $K$ with $j$-invariant $j$.
Let $p \ge 5$ and let $\fq \nmid p$ be a prime of $K$. 
Then $p \mid \# \overline{\rho}_{E,p}(I_\fq)$ if and only if $E$ has potentially multiplicative
reduction at $\fq$ (i.e.\ $\ord_\fq(j)<0$) and $p \nmid \ord_\fq(j)$.
\end{lem}
\begin{proof}
If $E$ has potentially good reduction at $\fq$ then 
$\# \overline{\rho}_{E,p}(I_\fq)$ divides $24$ (e.g.\ 
\cite[Introduction]{kraus2}). 
For $E$ with potentially multiplicative reduction
the result is a well-known consequence of the theory of the Tate
curve
(e.g.\ \cite[Proposition V.6.1]{SilvermanII}).
\end{proof}

\begin{lem}\label{lem:inertiaH}
Let $\fq \notin S$. 
Let $(a,b,c)$ be a solution to the Fermat equation \eqref{eqn:Fermat}
with prime exponent $p\ge 5$ such that $\fq \nmid p$. Let $E=E_{a,b,c}$
be the Frey curve in \eqref{eqn:Frey}. Then
$p \nmid \# \overline{\rho}_{E,p}(I_\fq)$. 
\end{lem}
\begin{proof}
By \eqref{eqn:inv}, if $a$, $b$, $c$ have unequal valuations at $\fq$, then
$\ord_\fq(j)<0$ and $p \mid \ord_\fq(j)$. Otherwise
$\ord_\fq(j) \ge 0$. In either case the lemma follows 
from Lemma~\ref{lem:inertiaGeneral}.
\end{proof}

\begin{lem}\label{lem:inertia3}
Let $E$ be an elliptic curve over $K$, and let $p \ge 3$.
Let $\mP \in S$ and suppose $E$ has potentially good reduction
at $\mP$. Let $\Delta$ be the discriminant of $E$ (not
necessarily minimal at $\mP$). Then $3 \mid \rho_{E,p}(I_\mP)$
if and only if $3 \nmid \ord_\mP(\Delta)$.
\end{lem}
\begin{proof}
Suppose $3 \nmid \ord_\mP(\Delta)$.
Let $L$ be the maximal
unramified extension of $K_\mP$. 
The valuation of $\Delta$ in $L$ is unchanged and so not divisible
by $3$. 
Now $E$ acquires good reduction over $L(E[p])$ (e.g.\ 
\cite[Proposition IV.10.3]{SilvermanII}).
It follows that the valuation of $\Delta$ in $L(E[p])$ is divisible by $12$.
Thus $3$ divides the degree $[L(E[p]):L]$.
However the Galois group of $L(E[p])/L$ is isomorphic to $\overline{\rho}_{E,p}(I_\mP)$, 
and so $3 \mid \# \overline{\rho}_{E,p}(I_\mP)$.

Kraus 
\cite[Th\'{e}or\`eme 3]{kraus2} 
enumerates all the 
possibilities for the Galois group of $L(E[p])/L$. 
If $3 \mid \ord_\mP(\Delta)$, its order 
is $1$, $2$, $4$ or $8$.
This proves the converse.
\end{proof}

\begin{lem}\label{lem:inertia}
Let $\mP \in S$. 
Let $(a,b,c)$ be a solution to the Fermat equation \eqref{eqn:Fermat}
with prime exponent $p> 4\ord_{\mP}(2)$. Let $E=E_{a,b,c}$
be the Frey curve in \eqref{eqn:Frey}. 
\begin{enumerate}
\item[(i)] If $\mP \in T$ then $E$ has potentially multiplicative reduction at $\mP$,
and
$p \mid \# \overline{\rho}_{E,p}(I_\mP)$. 
\item[(ii)] If $\mP \in U$ then
either $E$ has potentially multiplicative reduction at $\mP$ and
$p \mid \# \overline{\rho}_{E,p}(I_\mP)$, or 
$E$ has potentially good reduction at $\mP$ and $3 \mid \# \overline{\rho}_{E,p}(I_\mP)$. 
\end{enumerate}
\end{lem}
\begin{proof}
Let 
$\pi$ be a uniformizer
for $K_\mP$. 
Let
\[
t=\min\{ \ord_\mP(a), \ord_\mP(b), \ord_\mP(c) \}, \qquad
\alpha= \pi^{-t} a, \qquad \beta=\pi^{-t} b, \qquad \gamma=\pi^{-t} c.
\]
Then $\alpha$, $\beta$, $\gamma \in \OO_\pi$.
Suppose first that $\mP \in T$. By the definition of $T$,
the prime $\mP$ has residue field $\F_2$. 
As $\alpha^p+\beta^p+\gamma^p=0$,
precisely 
one of $\alpha$, $\beta$, $\gamma$ is divisible
by $\pi$. 
Thus $\ord_\mP(a)$, $\ord_\mP(b)$, $\ord_\mP(c)$ are not all equal;
two out of $a$, $b$, $c$ have valuation 
$t$ (say) and one has valuation $t+k$ with $k\ge 1$. 
From the formulae in \eqref{eqn:inv} we have
$\ord_\mP(j)=8 \ord_\mP(2)-2kp$. As $p> 4\ord_\mP(2)$, we see that 
$\ord_\mP(j)<0$ and $p \nmid \ord_\mP(j)$. Thus (i)
follows from Lemma~\ref{lem:inertiaGeneral}.

Suppose now that $\mP \in U$. 
If $\ord_\mP(a)$, $\ord_\mP(b)$, $\ord_\mP(c)$ are not all equal
then (ii) follows
as above. Suppose they are all equal. 
By the formulae in \eqref{eqn:inv} we have 
\[
\ord_\mP(j) \ge 8 \ord_\mP(2) > 0,
\qquad
\ord_\mP(\Delta)=4\ord_\mP(2)+6tp.
\] 
In particular,
$E$ has potentially good reduction. 
By definition of $U$, we have $3 \nmid \ord_\mP(2)$
and so $3 \nmid \ord_\mP(\Delta)$. 
Now (ii) follows from Lemma~\ref{lem:inertia3}.
\end{proof}

\section{Level Lowering and Eichler--Shimura}

\begin{thm}\label{thm:ll2}
Let $K$ be a totally real field satisfying (ES).
There is a constant $B_K$
depending only on $K$ such that 
the following hold.
Let $(a,b,c)$
be a non-trivial solution to the Fermat equation \eqref{eqn:Fermat}
with prime exponent $p>B_K$,
and rescale $(a,b,c)$ so that it remains integral and
satisfies $\cG_{a,b,c}=\fp$ for some $\fp \in \cH$.
Write $E$ for the Frey curve \eqref{eqn:Frey}.
Then there is an elliptic curve $E^\prime$ over $K$ such that
\begin{enumerate}
\item[(i)] the conductor of $E^\prime$ is divisible only by primes in 
$S \cup \{\fp\}$;
\item[(ii)] $\# E^\prime(K)[2]=4$;
\item[(iii)] $\overline{\rho}_{E,p} \sim \overline{\rho}_{E^\prime,p}$;
\end{enumerate}
Write $j^\prime$ for the $j$-invariant of $E^\prime$.
Then,
\begin{enumerate}
\item[(a)] for $\mP \in T$, we have $\ord_\mP(j^\prime)<0$;
\item[(b)] for $\mP \in U$, we have either $\ord_\mP(j^\prime)<0$
or $3 \nmid \ord_\mP(j^\prime)$;
\item[(c)] for $\fq \notin S$, we have $\ord_\fq(j^\prime) \ge 0$.
\end{enumerate}
In particular, $E^\prime$ has potentially good reduction away from $S$.
\end{thm}
\begin{proof}
We first observe, by Lemma~\ref{lem:cond}, that
$E$ is semistable outside $S \cup \{\fp\}$. 
By taking $B_K$ to be sufficiently large, we see
from Corollary~\ref{cor:Freymod} that $E$ is modular,
and from Theorem~\ref{thm:irred} that $\overline{\rho}_{E,p}$
is irreducible. 
Applying Theorem~\ref{thm:levell}
and Lemma~\ref{lem:cond} we see that $\overline{\rho}_{E,p} \sim \overline{\rho}_{\ff,\varpi}$
for a Hilbert newform $\ff$ of level $\cN_p$
and some prime $\varpi \mid p$ of $\Q_\ff$. Here $\Q_\ff$ is the field
generated by the Hecke eigenvalues of $\ff$. 

 Next we reduce to the case where 
$\Q_\ff=\Q$, after possibly enlarging $B_K$ by an effective amount.
This step uses standard ideas, originally due to Mazur
as indicated in the introduction, that can be found in 
\cite[Section 4]{BenS}, 
\cite[Proposition 15.4.2]{Cohen}, 
\cite[Section 3]{KrausFermat},
and so we omit the details.
We have assumed (ES): $T \ne \emptyset$, or $[K:\Q]$
is odd, or Conjecture~\ref{conj:ES} holds.
We would like to show there is some 
elliptic curve $E^\prime/K$ having the same $\mathrm{L}$-function as $\ff$.
This is immediate if we assume Conjecture~\ref{conj:ES}, and follows
from Theorem~\ref{thm:Zhang} if $[K:\Q]$ is odd. 
Suppose $T \ne \emptyset$ and let $\mP \in T$.
By Lemma~\ref{lem:inertia}, $E$ has 
potentially multiplicative reduction at $\mP$
and $p \mid \# \overline{\rho}_{E,p}(I_\mP)$. 
The existence of $E^\prime$ follows from Corollary~\ref{cor:ES}
after possibly enlarging $B_K$ to ensure $p \nmid (\norm_{K/\Q}(\mP) \pm 1)$.

We now know that $\overline{\rho}_{E,p} \sim \overline{\rho}_{E^\prime,p}$
for some $E^\prime/K$ with conductor $\cN_p$ given
by \eqref{eqn:cnp}.
After enlarging $B_K$ by an effective amount, 
and possibly replacing $E^\prime$ by an isogenous curve,
we may assume that $E^\prime$ has full $2$-torsion; this 
step 
again uses standard ideas 
\cite[Proposition 15.4.2]{Cohen}, 
\cite[Section 3]{KrausFermat},
and so we omit the details.

It remains to prove (a), (b), (c). There are finitely many 
elliptic curves $E^\prime$ with (full $2$-torsion and) good reduction
outside $S \cup \{\fp\}$. We may therefore, after possibly enlarging $B_K$,
suppose 
that for all primes $\fq$, if $\ord_\fq(j^\prime)  < 0$ then
$p \nmid \ord_\fq(j^\prime)$.
 Now we know from Lemma~\ref{lem:inertiaH}, for $\fq \notin S$,
that $p \nmid \# \overline{\rho}_{E,p}(I_\fq)$ and so 
$p \nmid \# \overline{\rho}_{E^\prime,p}(I_\fq)$. By Lemma~\ref{lem:inertiaGeneral}, we see that $\ord_\fq(j^\prime) \ge 0$, which proves (c).

For (a), let $\mP \in T$. Applying Lemma~\ref{lem:inertia}
we have $p \mid \overline{\rho}_{E,p}(I_\mP)$ and so 
$p \mid \overline{\rho}_{E^\prime,p}(I_\mP)$. 
Now by Lemma~\ref{lem:inertiaGeneral} we have $\ord_\mP(j^\prime)<0$
which proves (a).

For (b) suppose $\mP \in U$. 
If $p \mid \overline{\rho}_{E,p}(I_\mP)$
then again $\ord_\mP(j^\prime)<0$ as required. Thus 
suppose $p \nmid \overline{\rho}_{E,p}(I_\mP)$. By
Lemma~\ref{lem:inertia} we have $3 \mid \overline{\rho}_{E,p}(I_\mP)$. 
It follows that $p \nmid \overline{\rho}_{E^\prime,p}(I_\mP)$
and $3 \mid \overline{\rho}_{E^\prime,p}(I_\mP)$.
The first conclusion together with Lemma~\ref{lem:inertiaGeneral}
shows 
$\ord_\mP(j^\prime) \ge 0$ (recall, we have imposed
above the condition 
$\ord_\fq(j^\prime)<0$ implies $p \nmid \ord_\fq(j^\prime)$ for all $\fq$). 
By Lemma~\ref{lem:inertia3},
as $3 \mid \overline{\rho}_{E^\prime,p}(I_\mP)$,
we have $3 \nmid \ord_\mP(\Delta^\prime)$,
where $\Delta^\prime$ is the discriminant of 
$E^\prime$. But $j^\prime=(c_4^\prime)^3/\Delta^\prime$,
so $3 \nmid \ord_\mP(j^\prime)$ as required.
\end{proof}

\noindent \textbf{Remarks.}

\noindent \textbf{(I)} The constant $B_K$ is ineffective as it depends 
on the ineffective constant $A_K$ in Corollary~\ref{cor:Freymod}.
However, if we know that the Frey curve is modular (e.g.
if $K$ is real quadratic), the constant $B_K$
becomes effectively computable. Indeed, by the recipe in
\cite[Section 15.4]{Cohen}, all that is needed
are algorithms to compute 
the Hilbert newforms at the levels $\cN_p$,
the fields generated by their eigenvalues, 
and a finite number of eigenvalues for each,
and there are
effective algorithms for this \cite{DV}.

\medskip

\noindent \textbf{(II)} 
Theorem~\ref{thm:ll2}
is unconditional if $[K:\Q]$ is odd or $T \ne \emptyset$,
but is otherwise subject to Conjecture~\ref{conj:ES}.
So far as we know, Theorem~\ref{thm:Zhang}
is the strongest statement in the current literature
towards Conjecture~\ref{conj:ES}. It is fair to
ask if we can eliminate our dependence on Conjecture~\ref{conj:ES}
by using the full strength of Theorem~\ref{thm:Zhang}, and
the purpose of this remark is to explain why this seems
not possible, at least for quadratic fields.
Let $K=\Q(\sqrt{d})$
with $d>1$ squarefree
and $d \equiv 5 \pmod{8}$; these are the only real quadratic
fields with $T=\emptyset$.
We have in this case, through explicit calculations
that exploit the relation $\overline{\rho}_{E,p}\sim \overline{\rho}_{\ff,p}$, 
concocted
a mild set of local conditions on Fermat solution $(a,b,c)$,
at the primes in $S \cup \cH$, which if simultaneously
satisfied would ensure that $\pi_\fq$ is a principal
series at all finite primes $\fq$ of $K$, 
where $\pi$ is the automorphic representation attached to $\ff$.
We did not see a
way of showing that these local conditions cannot be
simultaneously satisfied.

\section{Proof of Theorem~\ref{thm:FermatGen}}
Theorem~\ref{thm:ll2} relates non-trivial solutions
of the Fermat equation (with $p$ sufficiently large) 
to elliptic curves $E^\prime$ with full $2$-torsion
having potentially good reduction outside $S$ and satisfying
certain additional properties. 
In this section we relate such
elliptic curves $E^\prime$ to solutions to the 
$S$-unit equation~\eqref{eqn:sunit}, using
basic facts about $\lambda$-invariants of elliptic curves
(e.g.\ \cite[Pages 53--55]{SilvermanI}), and use 
this to prove Theorem~\ref{thm:FermatGen}.
As $E^\prime$ has
full $2$-torsion over $K$, it has a model
\begin{equation}\label{eqn:e123}
E^\prime \; : \; y^2=(x-e_1)(x-e_2)(x-e_3).
\end{equation}
Here of course $e_1$, $e_2$, $e_3$ are distinct and so their
\textbf{cross ratio}
$\lambda=(e_3-e_1)/(e_2-e_1)$
belongs to $\PP^1(K)-\{0,1,\infty\}$.
Write $\sS_3$ for the symmetric group on $3$ letters.
The action of $\sS_3$
on $(e_1, e_2, e_3)$ extends via the cross ratio
to an action on $\PP^1(K)-\{0,1,\infty\}$, and
allows us to identify $\sS_3$ with the following subgroup of
$\PGL_2(K)$:
\[
\mathfrak{S}_3=\left\{z, \quad
1/z, \quad
1-z, \quad
1/(1-z),  \quad
z/(z-1), \quad
(z-1)/z
\right\} \, .
\]
The \textbf{$\lambda$-invariants of $E^\prime$} are the
the six elements (counted with multiplicity) of 
the $\sS_3$-orbit of $\lambda$; they 
are related to the $j$-invariant $j^\prime$ by 
\begin{equation}\label{eqn:j}
j^\prime=
2^8 \cdot (\lambda^2-\lambda+1)^3 \cdot \lambda^{-2} (\lambda-1)^{-2}.
\end{equation}

\begin{proof}[Proof of Theorem~\ref{thm:FermatGen}]
Let $K$ be a totally real field satisfying assumption (ES).
Let $B_K$ be as in Theorem~\ref{thm:ll2}, and let $(a,b,c)$
be a non-trivial solution to the Fermat equation~\eqref{eqn:Fermat}
with exponent $p>B_K$.
By Lemma~\ref{lem:gcd} we may rescale $(a,b,c)$ so that it remains
integral but $\cG_{a,b,c} =\fp$ for some $\fp \in \cH$. We now apply
Theorem~\ref{thm:ll2} which yields an elliptic curve $E^\prime/K$
with full $2$-torsion and potentially good reduction outside $S$
whose $j$-invariant $j^\prime$ satisfies the following two conditions:
\begin{enumerate}
\item[(a)] for all $\mP \in T$, we have $\ord_\mP(j^\prime) <0$;
\item[(b)] for all $\mP \in U$, we have $\ord_\mP(j^\prime)<0$
or $3 \nmid \ord_\mP(j^\prime)$.
\end{enumerate}
Let $\lambda$ be any of the $\lambda$-invariants of $E^\prime$.
Note that $j^\prime \in \OO_S$,
where $\OO_S$ is the ring of $S$-integers in $K$. By \eqref{eqn:j},
 $\lambda \in K$ satisfies a monic degree $6$ equation with
coefficients in $\OO_S$, and therefore 
$\lambda \in \OO_S$. However, $1/\lambda$, $\mu:=1-\lambda$ and $1/\mu$ are 
also a solutions to \eqref{eqn:j} and so belong to $\OO_S$.
Hence $(\lambda,\mu)$
is a solution to the $S$-unit equation \eqref{eqn:sunit}.
By assumption, this solution must satisfy either hypothesis
(A) or hypothesis (B) of Theorem~\ref{thm:FermatGen}. 
We now rewrite \eqref{eqn:j} as
\begin{equation}\label{eqn:jlammu}
j^\prime
=2^8 \cdot (1-\lambda \mu)^3 \cdot (\lambda \mu)^{-2} \, .
\end{equation}
and let
$t:=\max\{\lvert \ord_\mP(\lambda)\rvert, \lvert \ord_\mP(\mu)\rvert\}$.
Suppose that $(\lambda,\mu)$ 
satisfies (A): there is
some $\mP \in T$ so that 
$t \le 4 \ord_\mP(2)$.
If $t=0$ then $\ord_\mP(j^\prime) \ge 8 \ord_\mP(2) > 0$
which contradicts (a) above.
We may therefore suppose that $t>0$.
Now the relation $\lambda+\mu=1$
forces either $\ord_\mP(\lambda)=\ord_\mP(\mu)=-t$, or $\ord_\mP(\lambda)=0$ and $\ord_\mP(\mu)=t$, or $\ord_\mP(\lambda)=t$ and $\ord_\mP(\mu)=0$.
Thus $\ord_\mP(\lambda \mu)=-2t<0$ or $\ord_\mP(\lambda \mu)=t>0$.
Either way, 
$\ord_\mP(j^\prime)=8 \ord_\mP(2)-2 t \ge 0$,
which again contradicts (a).

Thus $(\lambda,\mu)$ satisfies (B):
there is some $\mP \in U$ such that
$t \le 4 \ord_\mP(2)$ 
and $\ord_\mP(\lambda \mu) \equiv \ord_\mP(2)
\pmod{3}$. The former implies $\ord_\mP(j^\prime) \ge 0$
as above, and the latter together with \eqref{eqn:jlammu}
gives $3 \mid \ord_\mP(j^\prime)$. This contradicts (b) 
completing the proof.
\end{proof}


\section{Proofs of Theorems~\ref{thm:FermatQuad} and~\ref{thm:d5mod8}}
We would like to understand the solutions to
\eqref{eqn:sunit} for
real quadratic $K$. 
Let
\begin{equation}\label{eqn:LambdaS}
\Lambda_S=\{(\lambda,\mu)\quad : \quad \lambda+\mu=1, \qquad \lambda,\; \mu \in \OO_S^*\}.
\end{equation}
The following two lemmas
 are easy consequences of the definitions.
\begin{lem}
The action of $\sS_3$ on $\PP^1(K)-\{0,1,\infty\}$ induces
an action on $\Lambda_S$
given by
$(\lambda,\mu)^{\sigma} :=(\lambda^{\sigma},1-\lambda^\sigma)$
for $(\lambda,\mu) \in \Lambda_S$ and $\sigma \in \sS_3$.
\end{lem}
\begin{lem}\label{lem:equiv}
Let $(\lambda_1,\mu_1)$ and $(\lambda_2,\mu_2)$ be elements of
$\Lambda_S$ belonging
to the same $\sS_3$-orbit. Then
\begin{enumerate}
\item[(i)]
$\max\{\lvert \ord_\mP(\lambda_1) \rvert,
\lvert \ord_\mP(\mu_1) \rvert \} 
 \le 4 \ord_\mP(2)$
if and only if
$\max\{\lvert \ord_\mP(\lambda_2) \rvert,
\lvert \ord_\mP(\mu_2) \rvert \} 
 \le 4 \ord_\mP(2)$.
\item[(ii)]
$\ord_\mP(\lambda_1 \mu_1) \equiv \ord_{\mP}(2) \pmod{3}$ if and only if
$\ord_\mP(\lambda_2 \mu_2) \equiv \ord_{\mP}(2) \pmod{3}$.
\end{enumerate}
\end{lem}
 We denote 
the set
of $\sS_3$-orbits in $\Lambda_S$
by $\sS_3 \backslash \Lambda_S$.
The
three elements $(2,-1)$, $(-1,2)$, $(1/2,1/2)$ of $\Lambda_S$ form a single
orbit. As stated in the introduction, we call these 
the \textbf{irrelevant solutions} to \eqref{eqn:sunit},
and call their orbit the 
\textbf{irrelevant orbit}.
Other solutions are called \textbf{relevant}.

In this section 
$K=\Q(\sqrt{d})$, where
$d\ge 2$ is a squarefree integer.
Note $T \ne \emptyset$ or $U \ne \emptyset$.
The 
irrelevant orbit satisfies
condition (A) of Theorem~\ref{thm:FermatGen}
if $T \neq \emptyset$, and satisfies
condition (B) if $U \ne \emptyset$. 

\begin{lem}\label{lem:makeint}
Let $K=\Q(\sqrt{d})$ where $d \ge 2$ is squarefree. 
Let $(\lambda,\mu) \in \Lambda_S$. 
\begin{enumerate}
\item[(i)] $\lambda$, $\mu \in \Q$ if and only if $(\lambda,\mu)$
belongs to the irrelevant orbit.
\item[(ii)] There is an element $\sigma \in \sS_3$
so that $(\lambda^\prime,\mu^\prime)=(\lambda,\mu)^\sigma$
satisfies $\lambda^\prime$, $\mu^\prime \in \OO_K$.
\end{enumerate}
\end{lem}
\begin{proof}
Note that $\lambda$ and $\mu$ are in $\Q$
if and only if they both have the form $\pm 2^r$.
From the relation $\lambda+\mu=1$ we quickly deduce (i).
For part (ii), 
as $K$ is quadratic, $\lvert S \rvert=1$ or $2$,
and the lemma follows from examining the possible
signs for $\ord_\mP(\lambda)$ and $\ord_\mP(\lambda^\sigma)$
for $\mP \in S$ and $\sigma \in \sS_3$.
\end{proof}

\begin{lem}\label{lem:param}
Let $K=\Q(\sqrt{d})$ where $d \ge 2$ is squarefree. 
Up to the action of $\sS_3$, every relevant $(\lambda,\mu) \in \Lambda_S$
has the form 
\begin{equation}\label{eqn:parasol}
\lambda=\frac{\eta_1 \cdot 2^{r_1}-\eta_2 \cdot 2^{r_2}+1 +v \sqrt{d}}{2},
\qquad 
\mu=\frac{\eta_2 \cdot 2^{r_2}-\eta_1 \cdot 2^{r_1}+1 -v \sqrt{d}}{2}
\end{equation}
where
\begin{equation}\label{eqn:paracond1}
\eta_1=\pm 1, \qquad \eta_2=\pm 1, \qquad r_1 \ge r_2 \ge 0, \qquad v \in \Z, \qquad v \ne 0
\end{equation}
are related by
\begin{gather}\label{eqn:paracond2}
(\eta_1 \cdot 2^{r_1}-\eta_2 \cdot 2^{r_2}+1)^2 
-\eta_1 \cdot 2^{r_1+2}
=
d v^2 ,
\\
\label{eqn:paracond3}
(\eta_2 \cdot 2^{r_2}-\eta_1 \cdot 2^{r_1}+1)^2 
-\eta_2 \cdot 2^{r_2+2}
=d v^2 .
\end{gather}
Moreover, if $d \not \equiv 1 \pmod{8}$ then we can take $r_2=0$.
\end{lem}
Observe that \eqref{eqn:paracond2} and \eqref{eqn:paracond3} are
equivalent, in the sense that we may rearrange either of the equations to 
obtain the other, but it is convenient to have both. 
\begin{proof}
Suppose $\eta_1$, $\eta_2$, $r_1$, $r_2$ and $v$ satisfy \eqref{eqn:paracond1}, \eqref{eqn:paracond2},
\eqref{eqn:paracond3}
and let $\lambda$, $\mu$ be given by \eqref{eqn:parasol}. 
It is clear that $\lambda$, $\mu$ belong to $\OO_S$ but not $\Q^*$,
and $\lambda+\mu=1$.
Moreover, from \eqref{eqn:paracond2}, \eqref{eqn:paracond3},
the norms of $\lambda$ and $\mu$ are $\eta_i 2^{r_i}$, and
thus $\lambda$, $\mu \in \OO_S^*$.
Hence $(\lambda,\mu)$
is a relevant element of $\Lambda_S$.
Conversely, suppose $(\lambda,\mu)$ is a relevant element of $\Lambda_S$. 
Thus $\lambda$, $\mu \notin \Q$, and
by Lemma~\ref{lem:makeint} we may suppose that $\lambda$, $\mu \in \OO_K$.
Let $x \mapsto \overline{x}$ denote conjugation in $K$.
Then
\[
\lambda \overline{\lambda}=\eta_1 \cdot 2^{r_1}, \qquad \mu \overline{\mu}=\eta_2 \cdot 2^{r_2}, \qquad
\eta_1=\pm 1, \qquad \eta_2=\pm 1.
\]
Swapping $\lambda$, $\mu$ if necessary (which does not change the orbit), 
we may suppose that $r_1 \ge r_2 \ge 0$.
Note that if $d \not \equiv 1 \pmod{8}$, then $S$ consists of one element,
and the  relation $\lambda+\mu=1$ forces $r_2=0$.
Now, 
\begin{equation*}
\lambda+\overline{\lambda} = \lambda \overline{\lambda} - (1-\lambda)(1-\overline{\lambda}) +1
						 = \lambda \overline{\lambda} - \mu \overline{\mu} +1 
						 = \eta_1 \cdot 2^{r_1}- \eta_2 \cdot 2^{r_2} +1 \, .
\end{equation*}
Moreover we can write $\lambda-\overline{\lambda}=v \sqrt{d}$, where $v \in \Z$, and as $\lambda \notin \Q$,
we have $v \ne 0$. The expressions for $\lambda+\overline{\lambda}$ and $\lambda-\overline{\lambda}$
give the expression for $\lambda$ in \eqref{eqn:parasol}, and we deduce the expression for $\mu$
from $\mu=1-\lambda$. Finally,
the identity
$(\lambda+\overline{\lambda})^2-(\lambda-\overline{\lambda})^2=4 \lambda \overline{\lambda}$ gives
 \eqref{eqn:paracond2} 
and the corresponding identity for $\mu$ gives
\eqref{eqn:paracond3}.
\end{proof}

\begin{table}
\begin{tabular}{||c|c|c||}
\hline \hline
$d$ & relevant elements of $\Lambda_S$  up to &        extra conditions \\
    & the action of $\sS_3$ and Galois conjugation  &  \\
\hline \hline
\multirow{2}{*}{$d=2$} &
 $(\sqrt{2},1-\sqrt{2})$, 
$(-16+12\sqrt{2},17-12 \sqrt{2})$,  & \\ 
& 
$(4+2\sqrt{2},-3+2\sqrt{2})$, $(-2+2\sqrt{2},3-2\sqrt{2})$ & \\
\hline
$d=3$ & 
$(2+\sqrt{3},-1-\sqrt{3})$, $(8+4\sqrt{3},-7-4\sqrt{3})$
& \\
\hline
\multirow{2}{*}{$d=5$}  & $\big( (1+\sqrt{5})/2,(1-\sqrt{5})/2\big)$, 
$(-8+4 \sqrt{5},9-4\sqrt{5})$, & \\ 
& $(-1+\sqrt{5},2-\sqrt{5})$ & \\
\hline
$d=6$ & $(-4+2\sqrt{6}, 5-2 \sqrt{6})$ & \\
\hline
$d \equiv 3 \pmod{8}$  & \multirow{2}{*}{none}  & \\
$d \ne 3$ & & \\
\hline
$d \equiv 5 \pmod{8}$  & \multirow{2}{*}{none} & \\
$d \ne 5$ & & \\
\hline
\multirow{2}{*}{$d \equiv 7 \pmod{8}$} & 
\multirow{2}{*}{
$( 2^{2s+1}+2^{s+1} w \sqrt{d},\, 1- 2^{2s+1}-2^{s+1}w \sqrt{d}   )$
}
 & $4^s-1=d w^2$\\
  & & $s \ge 2$, $w \ne 0$ \\
\hline
$d \equiv 2 \pmod{16}$ &
\multirow{2}{*}{
$
( - 2^{2s}+2^s w \sqrt{d},
1+2^{2s}- 2^s w \sqrt{d}
)
$
}
 & $4^s+2=d w^2$ \\
$d \ne 2$ & & $s \ge 2$, $w \ne 0$\\
\hline
$d \equiv 6 \pmod{16}$ & \multirow{2}{*}{none} & \\
$d \ne 6$ & & \\
\hline
$d \equiv 10 \pmod{16}$ & none & \\
\hline
\multirow{2}{*}{$d \equiv 14 \pmod{16}$} & 
\multirow{2}{*}{
$
( 2^{2s}+2^s w \sqrt{d},
1-2^{2s}- 2^s w \sqrt{d}
)
$
}

& $4^s-2=d w^2$\\
   & & $s \ge 2$, $w \ne 0$\\
\hline \hline
\end{tabular}
\caption{The relevant elements of $\Lambda_S$ for 
$d \ge 2$ squarefree, $d \not \equiv 1 \pmod{8}$.}
\label{table:sols}
\end{table}

\begin{lem}\label{lem:table}
Let $d \not \equiv 1 \pmod{8}$ be squarefree $\ge 2$.
The relevant elements of $\Lambda_S$, up to the 
action of $\sS_3$ and Galois conjugation are as given in
Table~\ref{table:sols}.
\end{lem}
\begin{proof}
We apply Lemma~\ref{lem:param}. 
As $d \not \equiv 1 \pmod{8}$,
we have $r_2=0$. Values $0 \le r_1 \le 5$ (together with $\eta_1=\pm 1$, $\eta_2=\pm 1$)
in \eqref{eqn:paracond2} yield the solutions 
given in Table~\ref{table:sols} for $d=2$, $3$, $5$ and $6$,
and also the solution $(16-4 \sqrt{14},-15+4\sqrt{14})$ which
is included under the table entry for $d \equiv 14 \pmod{16}$.
We may therefore suppose that $r_1 \ge 6$. 
If $\eta_2=-1$ then \eqref{eqn:paracond3} gives
$2^{2 r_1} +4= d v^2$. As $d$ is squarefree and $r_1 \ge 6$ this gives 
$d \equiv 1 \pmod{8}$,
which is a contradiction. Thus $\eta_2=1$. 
Now \eqref{eqn:paracond2} gives
\[
2^{r_1+2}(2^{r_1-2} - \eta_1)  = d v^2.
\]
As $d$ is squarefree, 
we see that $r_1$ is even precisely when $d$ is odd.
Suppose first that $d$ is odd, and so
$r_1=2s+2$ and $v=2^{s+2} w$ 
for some non-zero integer $w$. 
As $r_1 \ge 6$, we have $s \ge 2$. 
Then $4^s-\eta_1=d w^2$. As $\eta_1 = \pm 1$,
this equation is impossible
if $d \equiv 3$ or $5 \pmod{8}$. This completes
the proofs of the entries in the table for $d \equiv 3$, $5 \pmod{8}$
(including $d=3$, $d=5$).
If $d \equiv 7 \pmod{8}$, then reducing mod $8$ shows 
that $\eta_1=1$. This gives the solution in Table~\ref{table:sols}
for $d \equiv 7 \pmod{8}$.

Suppose now that $d$ is even, and so $r_1=2s+1$ and $v=2^{s+1} w$ 
for some non-zero integer $w$, and $s \ge 2$.
Then $2^{2s-1}-\eta_1 = (d/2) w^2$. This gives a contradiction modulo $8$,
if $d/2 \equiv 3$ or $5 \pmod{8}$, and proves the correctness of the
entries in the table for $d \equiv 6$, $10 \pmod{16}$ (including $d=6$).
Moreover, if $d \equiv 2$ or $14 \pmod{16}$ then $\eta = -1$, $1$
respectively, completing the proof except for $d=2$.

For $d=2$, we have $2^{2s-1}+1=w^2$. The factorization $(w-1)(w+1)=2^{2s-1}$
quickly leads to a contradiction, so there are no further solutions. 
\end{proof}

\subsection{Proof of Theorem~\ref{thm:FermatQuad}}
For now let $d \ge 2$ be squarefree and satisfy $d \not \equiv 1 \pmod{4}$.
Write $K=\Q(\sqrt{d})$. 
Then $S=T=\{ \mP \}$ say; in particular assumption (ES) is satisfied.
We apply Theorem~\ref{thm:FermatGen} and focus on
verifying condition (A). 
In view of Lemma~\ref{lem:equiv}, we only have to
verify (A) for one representative of each
$\sS_3$-orbit. It is easy to check that the irrelevant orbit,
as well as the solutions listed in Table~\ref{table:sols}
for $d=2$, $3$, $6$ all satisfy (A).
For $d \equiv 3 \pmod{8}$ or $d \equiv 6$, $10 \pmod{16}$
there are no further solutions and so the proof is complete
for parts (i), (ii) of Theorem~\ref{thm:FermatQuad}.
Suppose $d \equiv 2 \pmod{16}$, and $d \ne 2$. 
We know from Table~\ref{table:sols}
that the $S$-unit equation has no relevant solutions,
unless there are $s \ge 2$ and $w \ne 0$ such that $4^s+2=d w^2$.
Now if $q \mid d$ is an odd prime, then $-2$ is
a quadratic residue modulo $q$ and so $q \equiv 1$, $3 \pmod{8}$.
This proves (iii)
of Theorem~\ref{thm:FermatQuad}, and (iv)
is similar. 

\subsection{Proof of Theorem~\ref{thm:d5mod8}}
Here $d \equiv 5 \pmod{8}$. 
Again we apply Theorem~\ref{thm:FermatGen}, 
but check condition (B) of this theorem. Note $U=S=\{\mP\}$
where $\mP=2\cdot \OO_K$. Moreover, the irrelevant solutions to the
$S$-unit equation~\eqref{eqn:sunit} satisfy condition (B).
From Table~\eqref{table:sols}, there are no further
solutions to the $S$-unit equation for $d \ne 5$. This
completes the proof. 
For $d=5$, Table~\eqref{table:sols} lists three 
solutions to the $S$-unit equation 
that satisfy $\ord_\mP(\lambda \mu)=0$, $2$, $1$ respectively;
we cannot complete the proof in this case, which is why
it is excluded from the statement of the theorem. 

\section{Proof of Theorem~\ref{thm:density}}\label{sec:density}
For a set $\cU$ of positive integers
and a positive real number $X$, we let $\cU(X)=\{d \in \cU : d \le X\}$.
We define the \textbf{absolute density} of $\cU$ 
to be
\[
\delta(\cU)=\lim_{X \rightarrow \infty}
\# \cU(X)/X,
\]
provided the limit exists.
In this section, as in the Introduction, we let $\Nsf$ be the set
of squarefree integers $d \ge 2$. For $\cU \subseteq \Nsf$,
recall the definition of the \textbf{relative density} of $\cU$ in $\Nsf$
given in \eqref{eqn:deltarel}, which we can now write as
\[
\deltarel(\cU)=\lim_{X \rightarrow \infty}
  \# \cU(X) /  \#  \Nsf (X)
\]
provided the limit exists. 
We need the following two classical analytic theorems.
\begin{thm} \textup{(e.g.\ \cite[page 636]{Landau})} \label{thm:sfden}
For integers $r$, $N$ with $N$ positive, let 
\[
\Nsf_{r,N}=\{ d \in \Nsf \; : \; d \equiv r \pmod{N}\}.
\]
Let $s=\gcd(r,N)$ and suppose that $s$ is squarefree.
Then 
\[
\# \Nsf_{r,N}(X)
\sim
\frac{\varphi(N)}{s \varphi(N/s) N \prod_{q \mid N} (1-q^{-2})} \cdot \frac{6}{\pi^2} X \, .
\]
\end{thm}
Here $\varphi$ denotes Euler's totient function.

\begin{thm} \textup{(e.g.\ \cite[pages 641--643]{Landau})} \label{thm:ana2}
Let $N$ be a positive integer. Let $r_1,\dotsc, r_m$
be distinct modulo $N$, satisfying $\gcd(r_i,N)=1$. 
Let $E$ be the set of positive integers $d$
such that every prime factor $q$ of $d$ satisfies
$q \equiv r_i \pmod{N}$ for some $i=1,\dotsc,k$.
Then there is some positive constant $\gamma=\gamma(N,r_1,\dots,r_m)$ such that
\[
\# E(X) \sim  \gamma \cdot X / \log(X)^{1-\frac{m}{\varphi(N)}} \, .
\]
\end{thm}

Let $\sC$ and $\cD$ given in \eqref{eqn:cD}.
\begin{lem}\label{lem:crucial}
Let
$\sC^\prime= \Nsf \backslash \sC$. Then
$\delta(\sC^\prime)=0$.
\end{lem}
The proof of Lemma~\ref{lem:crucial} is somewhat lengthy.
Before embarking on the proof, we show that the lemma
is enough to imply Theorem~\ref{thm:density}.

\begin{proof}[Proof of Theorem~\ref{thm:density}]
Observe that Theorem~\ref{thm:density} follows immediately from
Theorem~\ref{thm:FermatGen} provided we can
prove the density claims in \eqref{eqn:density}.
Theorem~\ref{thm:sfden} applied to $\Nsf=\Nsf_{0,1}$ gives 
$\# \Nsf(X) \sim  6 X/\pi^2$.
It follows for $\cU \subseteq \Nsf$ that $\delta(\cU)$
exists if and only if $\deltarel(\cU)$ exists, and in this case, the two
are related by
$\deltarel(\cU)=\pi^2 \delta(\cU)/6$.
By Lemma~\ref{lem:crucial} we have  
$\deltarel(\sC^\prime)=\delta(\sC^\prime)=0$. 
As $\sC^\prime$, $\sC$ are complements in $\Nsf$,
we have $\deltarel(\sC)=1$. 
It remains to prove $\deltarel(\cD)=5/6$. 
By definition, $\cD=\sC \cap (\Nsf- \Nsf_{5,8})$, so it suffices to prove 
$\deltarel(\Nsf_{5,8})=1/6$. This follows from Theorem~\ref{thm:sfden}.
\end{proof}

\subsection{Proof of Lemma~\ref{lem:crucial}}
To complete the proof of Theorem~\ref{thm:density} we need to 
prove Lemma~\ref{lem:crucial}.
By definition, $\sC^\prime$ is the set of squarefree $d \ge 2$ such that
the $S$-unit equation~\eqref{eqn:sunit} has a 
relevant solution in $\Q(\sqrt{d})$. 
By Lemma~\ref{lem:param} this is precisely the set
of squarefree $d \ge 2$ satisfying equations~\eqref{eqn:paracond2}
and~\eqref{eqn:paracond3} where $\eta_1=\pm 1$, $\eta_2=\pm 1$,
$r_1 \ge r_2 \ge 0$ and $v \ne 0$. For $\eta_1$, $\eta_2=\pm 1$, 
and $\kappa_1$, $\kappa_2=0$, $1$ write
$\sC^\prime(\eta_1,\eta_2,\kappa_1,\kappa_2)$ for the set
of $d \ge 2$ satisfying equations~\eqref{eqn:paracond2}
and~\eqref{eqn:paracond3} with $r_i \equiv \kappa_i \pmod{2}$.
Then $\sC^\prime$ is the union of all sixteen
$\sC^\prime({\eta_1,\eta_2,\kappa_1, \kappa_2})$.
Let $q$ be an odd prime divisor of $d$. 
If $d \in \sC^\prime(\eta_1,\eta_2,\kappa_1,\kappa_2)$ then
reducing equations~\eqref{eqn:paracond2}
and~\eqref{eqn:paracond3} modulo $q$ shows that
$\eta_1 2^{\kappa_1}$ and $\eta_2  2^{\kappa_2}$
are both quadratic residues modulo $q$.
If $(\eta_1,\eta_2,\kappa_1,\kappa_2) \ne (1,1,0,0)$,
then the possible odd prime divisors $q$ of $d$
belong to a {\em proper} subset of the congruence classes
$\{\overline{1},\overline{3},\overline{5},\overline{7}\}$
modulo $8$. 
It follows from Theorem~\ref{thm:ana2} that
$\delta(\sC^\prime (\eta_1,\eta_2,\kappa_1,\kappa_2))=0$ for
$(\eta_1,\eta_2,\kappa_1,\kappa_2) \ne (1,1,0,0)$.

To show $\delta(\sC^\prime)=0$ it is now enough to show
that $\delta(\sC^\prime(1,1,0,0))=0$. Suppose $d \in \sC^\prime(1,1,0,0)$.
Thus there is a solution to \eqref{eqn:paracond2}, \eqref{eqn:paracond3}
with $\eta_1=\eta_2=1$, $r_1=2s$, $r_2=2t$ and $s \ge t \ge 0$.
Equations \eqref{eqn:paracond2}  and \eqref{eqn:paracond3}
are now equivalent to
\begin{equation}\label{eqn:fourprod}
(2^s+2^t+1)(2^s+2^t-1)(2^s-2^t+1)(2^s-2^t-1)=dv^2.
\end{equation}
Since $d>0$
we in fact have $s>t$. 
We write $\sC^\prime(1,1,0,0)=\cL \cup \cM$ where $\cL$
is the subset of $d$ for which there is a solution
to \eqref{eqn:fourprod} 
with
$t>0$, and $\cM$ is the subset of $d$
for which there is a solution $t=0$.
We separately show that $\delta(\cL)=0$ and $\delta(\cM)=0$.
Suppose now that $d \in \cL$.
Write $\mathbf{s}=(s,t)$,
and let
\begin{equation}\label{eqn:alphas}
\alpha_{1,\bs}=2^s+2^t+1, \quad \alpha_{2,\bs}=2^s+2^t-1,
\quad \alpha_{3,\bs}=2^s-2^t+1, \quad \alpha_{4,s}=2^s-2^t-1.
\end{equation}
These are odd and pairwise coprime.
Then $\alpha_{i,\bs}=d_{i,\bs} v_{i,\bs}^2$,
where $d_{i,\bs}$ are squarefree, and 
$d=\prod d_{i,\bs}$. Thus
\[
\cL=\{ d_{1,\bs} d_{2,\bs} d_{3,\bs} d_{4,\bs} \; : \; \bs=(s,t), \qquad s>t>0  \}.
\]
Recall we want to show that $\delta(\cL)=0$.
In fact we prove the  equivalent statement that $\deltasup(\cL)=0$,
where
$\deltasup(\cL)=\limsup_{X \rightarrow \infty} \, \# \cL(X)/X$.

For a positive integer $m$, we write $M_m=2^m-1$; 
this is the $m$-th Mersenne number. We  
make frequent use of the fact that $M_n \mid M_m$ whenever $n \mid m$.
\begin{lem}\label{lem:alphai}
Let $m$ be a positive integer. Let $M_m=2^m-1$.
Suppose $\bs_1 \equiv \bs_2 \pmod{m}$. Then $\alpha_{i,\bs_1} \equiv \alpha_{i,\bs_2} \pmod{M_m}$
for $i=1,\dots,4$.
\end{lem}
\begin{proof}
Write $\bs_1=(s_1,t_1)$, $\bs_2=(s_2,t_2)$. If $\bs_1 \equiv \bs_2 \pmod{m}$
then $f=\lvert s_1-s_2 \rvert$ and $g=\lvert t_1-t_2\rvert$ are divisible by $m$. But
$\alpha_{i,\bs_1}-\alpha_{i,\bs_2}$ can be written as a linear combination of $2^f-1$ and
$2^g-1$, and is hence divisible by $M_m$.
\end{proof}
\begin{lem}\label{lem:den1}
Let $m$ be a positive integer. Let $s_0>t_0>0$ and write $\bs_0=(s_0,t_0)$.
Let 
\[
\cA_{m,\bs_0}=
\{
d_{1,\bs} \; : \; \bs=(s,t), \quad s>t>0,
\quad \bs \equiv \bs_0 \pmod{m}
\} \, .
\]
Let $p_1,\dotsc,p_k$ be the distinct primes dividing $M_m$ that do not
divide $\alpha_{1,\bs_0}$, and write $N=p_1\cdots p_k$. 
Then, 
\[
\# \cA_{m,\bs_0} (X) \le 2^{-k} \cdot X+N.
\]
\end{lem}
\begin{proof}
Clearly $N \mid M_m$,
and $\alpha_{1,\bs_0}$ is coprime to $N$.
Suppose $\bs \equiv \bs_0 \pmod{m}$.
By Lemma~\ref{lem:alphai},
$\alpha_{1,\bs} \equiv \alpha_{1,\bs_0} \pmod{N}$, and so $\alpha_{1,\bs}$
is also coprime to $N$. As 
$\alpha_{i,\bs}=d_{i,\bs} v_{i,\bs}^2$, we have
$d_{1,\bs} \equiv \alpha_{1,\bs_0} v_{1,\bs}^{-2} \pmod{N}$.
Thus $d_{1,\bs}$ modulo $N$ belongs to the set 
$\{ \alpha_{1,\bs_0} \cdot w^2 : w \in (\Z/N\Z)^*\}$.
As $N$ is squarefree with $k$ distinct odd prime factors, 
this set
has cardinality $\varphi(N)/2^k$. 
The lemma follows. 
\end{proof}

We denote the number of distinct prime divisors
of a positive integer $n$ by $\omega(n)$.
\begin{lem}\label{lem:den2}
For $m \ge 1$, let $h_m=\omega(M_m)$. Then,
$\deltasup(\cL) \le {2^{-h_m/2}} \cdot m^2$.
\end{lem}
\begin{proof}
For now, fix some $\bs_0$ and let $k$ and $N$ be as in Lemma~\ref{lem:den1}. Let
\[
\cL_{m,\bs_0}=
\{
d_{1,\bs} d_{2,\bs} d_{3,\bs} d_{4,\bs} \; : \; \bs=(s,t), \quad s>t>0,
\quad \bs \equiv \bs_0 \pmod{m} 
\}.
\]
Then
\[
\# \cL_{m,\bs_0}(X) \le \# \cA_{m,\bs_0}(X)
\le 2^{-k}  {X}+N \le 2^{-k}  {X}+M_m,
\]
where the first inequality is clear from the definitions of $\cL_{m,\bs_0}$ and $\cA_{m,\bs_0}$,
and the second and third follow from Lemma~\ref{lem:den1}.
Now let $q_1,\dotsc,q_{k^\prime}$ be the distinct prime divisors of $M_m$ that do not divide $\alpha_{2,\bs_0}$.
Then (similarly to the above) we have
\[
\# \cL_{m,\bs_0}(X) 
\le 2^{-k^\prime} {X}+M_m.
\]
As $\alpha_{1,\bs_0}$ and $\alpha_{2,\bs_0}$ are coprime, either $k\ge h_m/2$ or $k^\prime\ge h_m/2$.
Hence
\[
\# \cL_{m,\bs_0}(X) 
\le {2^{-h_m/2}} X+M_m.
\]
Now let $\bs_1,\dotsc,\bs_{m^2}$ be a complete system of representatives for $\bs$ modulo $m$.
Then $\cL$ is the disjoint union of $\cL_{m,\bs_1},\dots,\cL_{m,\bs_{m^2}}$. 
Thus
\[
\# \cL(X) 
\le 
{2^{-h_m/2}}
\cdot
{m^2\cdot X}
+m^2 M_m.
\]
\end{proof}
A prime $\ell$ is a \textbf{primitive divisor}
of $M_m$ if $\ell \mid M_m$ but $\ell \nmid M_{m^\prime}$ for
all $m^\prime<m$. 
The following 
is a 
special case of the celebrated Primitive Divisor
Theorem of Bilu, Hanrot and Voutier \cite{BHV}.
\begin{thm}
If $n \ne 1$, $6$, then $M_n$ has a primitive divisor.
\end{thm}
\begin{cor}\label{cor:wm}
$h_m \ge 2^{\omega(m)}-2$.
\end{cor}
\begin{proof}
The number of divisors $n \mid m$ is at least $2^{\omega(m)}$.
For $n \ne 1$, $6$ dividing $m$,
let $\ell_n$ be a primitive
divisor of $M_n$. Then $\ell_n$ divides $M_m$,
and by definition $\ell_n \ne \ell_{n^\prime}$
whenever $n \ne n^\prime$.
This gives
at least $2^{\omega(m)}-2$ distinct primes dividing $M_m$.
\end{proof}

\begin{lem}\label{lem:cL}
$\delta(\cL)=0$.
\end{lem}
\begin{proof}
Combining
Lemma~\ref{lem:den2} and Corollary~\ref{cor:wm}, we have
\begin{equation}\label{eqn:lim}
0 \le \deltasup(\cL) \le
\frac{m^2}{2^{(h_m/2)}} \le \frac{m^2}{2^{(2^{\omega(m)-1} -1)}} 
\end{equation}
for any $m \ge 1$. Now let $y$ be large, and let
$m=\prod_{p \le y} p$. Then 
$\omega(m)=\pi(y)$ and
$m=\exp(\vartheta(y))$,
where
$\pi$ and $\vartheta$ are respectively the prime counting function
and the first Chebyshev function. 
By the Prime Number Theorem (e.g.\ \cite[Chapter 4]{Apostol}) 
\[
\pi(y) \sim {y}/{\log{y}}, \qquad \vartheta(y) \sim y.
\]
Letting $y \rightarrow \infty$ in \eqref{eqn:lim} 
clearly gives $\deltasup(\cL)=0$,
as required.
\end{proof}

Recall that we had written $\sC^\prime(1,1,0,0)=\cL \cup \cM$,
and wanted to show that $\delta(\cL)=\delta(\cM)=0$.
The following completes the proof of Lemma~\ref{lem:crucial}
and Theorem~\ref{thm:density}.
\begin{lem}\label{lem:cM}
$\delta(\cM)=0$.
\end{lem}
\begin{proof}
The set $\cM$ is the set of squarefree $d \ge 2$ such that \eqref{eqn:fourprod}
holds with $s>0$ and $t=0$. Then
$(2^{s-1}-1)(2^{s-1}+1)=d w^2$
where $w=v/2^{s+1} \in \Z\backslash \{0\}$. Now we can apply
a straightforward simplification of the above argument 
to show $\delta(\cM)=0$.
\end{proof}
We thank 
A.\ Meyerowitz,
G.\ Helms and other 
 users of {\tt mathoverflow.net} 
\footnote{http://mathoverflow.net/questions/149511/squarefree-parts-of-mersenne-numbers}
for suggesting alternative
proofs for Lemma~\ref{lem:cM}.
It is clear however that the ideas 
in these alternative proofs cannot be adapted to prove
Lemma~\ref{lem:cL}.

\end{document}